\documentclass[titlepage,11pt]{article}
\usepackage{latexsym}
\usepackage{amsfonts}
\usepackage{amsmath}
\usepackage{tikz}
\usepackage{float}
\usepackage{comment}
\newcommand\blackslug{\hbox{\hskip 1pt \vrule width 4pt height 8pt depth 1.5pt
        \hskip 1pt}}
\newcommand\bbox{\hfill \quad \blackslug \bigbreak}
\def\DD{\hbox{-}}
\def\CC{\hbox{-}\cdots\hbox{-}}
\def\LL{,\ldots,}
\def\cupcup{\cup\cdots\cup}

\title{Proof of a conjecture of Plummer and Zha}
\author{Maria Chudnovsky\thanks{Supported by NSF DMS-EPSRC grant DMS-2120644.}\\
Princeton University, Princeton, NJ 08544
\\
\\
Paul Seymour\thanks{Supported by AFOSR grant
A9550-19-1-0187.}\\
Princeton University, Princeton, NJ 08544}

\date{January 8, 2022; revised \today}

\newtheorem{thm}{}[section]

\newcommand{\Proof}{\noindent{\bf Proof.}\ \ }

\begin{document}
\maketitle
\begin{abstract}
Say a graph $G$ is a {\em pentagraph} if every cycle has length at least five, and every induced cycle of
odd length has length five. N. Robertson proposed the conjecture that the Petersen graph is the only pentagraph that
is three-connected and internally 4-connected, but this was disproved by M. Plummer and X. Zha in 2014.
Plummer and Zha conjectured that every  3-connected, internally 4-connected pentagraph is three-colourable.
We prove this: indeed, we will prove that every pentagraph is three-colourable.

\end{abstract}

\section{Introduction}
Let us say a graph $G$ is a {\em pentagraph} if every cycle has length at least five, and every induced cycle of
odd length has length five. (All graphs in this paper are finite, and have no loops or parallel edges.) 
Such graphs seem to be richly structured; indeed N. Robertson~\cite{neil} proposed the conjecture that the Petersen graph 
is the only pentagraph that
is three-connected and internally 4-connected, although this was disproved by M. Plummer and X. Zha~\cite{plummer} 
(see also~\cite{nelson}).

In the same paper, Plummer and Zha proposed the conjecture that all 3-connected, internally 4-connected pentagraphs have bounded
chromatic number, and the stronger conjecture that they are all three-colourable. The first was proved by Xu, Yu, and 
Zha~\cite{XYZ}, who proved that all pentagraphs are four-colourable; and we will prove the second.
Our main theorem is:
\begin{thm}\label{mainthm}
Every pentagraph is three-colourable.
\end{thm}

The {\em girth} of $G$ is the minimum length of a cycle in $G$; a {\em hole} in $G$ is an induced cycle of length at least 
four, and an
{\em odd hole} means a hole with odd length; and if $X\subseteq V(G)$, $G[X]$ denotes the subgraph induced
on $X$.
We first prove the result of~\cite{XYZ}, because the proof is short and pretty:
\begin{thm}\label{fourthm}
Every pentagraph is four-colourable.
\end{thm}
\Proof
It is enough to show that every connected pentagraph is four-colourable. Let $G$ be a connected pentagraph, let $v_0\in V(G)$,
and for each $k\ge 0$ let $L_k$ be the set of vertices with distance exactly $k$ from $v_0$. Thus the sets $L_0, L_1,\ldots$ 
are pairwise disjoint and have union $V(G)$; and for each $k\ge 1$,
each vertex in $L_k$ has a neighbour in $L_{k-1}$ and has no neighbour in $L_0\cupcup L_{k-2}$.
Suppose that for some $k\ge 0$, 
$G[L_k]$ 
is not bipartite, and choose a minimum such value of $k$. Thus $k\ge 1$; and since $G[L_k]$ is not bipartite, it contains 
an odd cycle
as a subgraph, and hence an induced odd cycle, and hence a hole $C$ of length five (since $G$ is
a pentagraph). Let $C$  have vertices $c_1\DD c_2\DD c_3\DD c_4\DD c_5\DD c_1$ in order. For $1\le i\le 5$,
let $d_i\in L_{k-1}$ be adjacent to $c_i$. It follows that each of $d_1\LL d_5$ has only one neighbour in
$V(C)$, since $G$ has girth at least five, and hence $d_1\LL d_5$ are all distinct.
Thus $k\ge 2$. Since $G[L_0\cupcup L_{k-2}]$ is connected, and $d_1, d_3$ 
both have a neighbour in $L_{k-2}$, there is an induced path $P$ between $d_1,d_3$ with interior in $L_0\cupcup L_{k-2}$.
Since 
$$d_1\DD P\DD d_3\DD c_3\DD c_4\DD c_5\DD c_1\DD d_1$$ is a hole of length at least six, it has even length and so 
$P$ has odd length. Consequently 
$$d_1\DD P\DD d_3\DD c_3\DD c_2\DD c_1\DD d_1$$ is an odd hole, and so it has
length five, and therefore $d_1d_3$ is an edge. Similarly $d_3d_5,d_5d_2,d_2d_4, d_4d_1$
are edges, and so $G[L_{k-1}]$ has a cycle of length five, contradicting the choice of $k$.

Thus $G[L_k]$ is bipartite for each $k$, and so $G$ is four-colourable. This proves \ref{fourthm}.~\bbox

Now we turn to the proof of \ref{mainthm}.
This is a consequence of a stronger statement, that every non-bipartite pentagraph is either isomorphic to the Petersen graph, 
or has a vertex of degree at most two,
or admits one of two kinds of decomposition, that a minimal non-three-colourable pentagraph cannot admit.
Let us see these decompositions. 

\begin{itemize}
\item A {\em $P_3$-cutset} means an induced three-vertex path $P$ of $G$ such that $G\setminus V(P)$ is disconnected.
\item A {\em parity star-cutset} means a set $X\subseteq V(G)$ such that $G\setminus X$ is disconnected, and there is a 
vertex $x\in X$ such that $x$ is adjacent
to every other vertex in $X$, and there is a component $A$ of 
$G\setminus X$ such that every two vertices in $X\setminus \{x\}$ are joined by an induced path of even length with 
interior in $V(A)$. We call this a {\em strong parity star-cutset} if $A$ can be chosen such that
in addition $x$ has a neighbour in $V(A)$.
\end{itemize}
We will prove:
\begin{thm}\label{decompthm}
Let $G$ be a pentagraph. Then either 
\begin{itemize}
\item $G$ is bipartite; or
\item $G$ is isomorphic to the Petersen graph; or
\item $G$ has a vertex of degree at most two; or
\item $G$ admits a $P_3$-cutset or a strong parity star-cutset.
\end{itemize}
\end{thm}
\noindent{\bf Proof of \ref{mainthm}, assuming \ref{decompthm}.\ \ }
We prove by induction on $|V(G)|$ that every pentagraph is three-colourable. Let $G$ be a pentagraph such that every
pentagraph with fewer vertices is three-colourable. If $G$ is isomorphic to the Petersen graph, or
some vertex has degree at most two, then $G$ is three-colourable; so by \ref{decompthm} we may assume that $G$ admits a
$P_3$-cutset or a parity star-cutset (indeed, a strong parity star-cutset, but we do not need ``strong'' here).

Suppose first that $G$ admits a $P_3$-cutset, and let $v_1\DD v_2\DD v_3$
be an induced path such that $G\setminus \{v_1,v_2,v_3\}$ is disconnected. Let $A_1$ be the union of at least one and not all
of the components of $G\setminus \{v_1,v_2,v_3\}$,
and let $A_2$ be the union of all the other components. For $i = 1,2$ let $G_i=G[A_i\cup \{v_1,v_2,v_3\}]$.
From the inductive hypothesis, both $G_1$ and $G_2$
are three-colourable; let $\phi_i:V(G_i)\rightarrow\{1,2,3\}$ be a three-colouring, for $i = 1,2$. We may assume that 
$\phi_i(v_1)=1$ and $\phi_i(v_2)=2$ for $i = 1,2$. Thus $\phi_1(v_3), \phi_2(v_3)\in \{1,3\}$, and if $\phi_1(v_3)= \phi_2(v_3)$
then $G$ is three-colourable. Thus we may assume that $\phi_1(v_3)=1$ and $\phi_2(v_3)=3$. Let $H_1$ be the subgraph of
$G_1$ induced on the set of vertices $v\in V(G_1)$ with $\phi_1(v)\in \{1,3\}$. If $v_1,v_3$ belong to different components 
of $H_1$,
then by exchanging colours in the component containing $v_3$, we obtain another three-colouring of $G_1$ that can be combined
with $\phi_2$ to show that $G$ is three-colourable. So we may assume that $v_1,v_3$ belong to the same component of $H_1$,
and so there is an induced path $P_1$ of $H_1$ between $v_1,v_3$. Consequently $P_1$ has even length, and length at least 
four since $G$ has girth at least five. Define $H_2$ in $G_2$ similarly: then similarly we may assume that $v_1,v_3$
belong to the same component of $H_2$, and so there is an induced path $P_2$ of $H_2$ between $v_1,v_3$ with odd length, at 
least three. But then $P_1\cup P_2$ is an induced cycle of $G$ of odd length at least seven, a contradiction.

Now suppose that $G$ admits a parity star-cutset, and let $X\subseteq V(G)$ and $v\in V(G)\setminus X$, such that 
$v$ is adjacent
to every vertex in $X$, and $G\setminus (X\cup \{v\})$ is disconnected, and there is a component $A$ of
$G\setminus (X\cup \{v\})$ such that every two vertices in $X$ are joined by an induced path of even length with 
interior in $V(A)$. Choose $X$ minimal with this property, and let $A_1\LL A_k$ be the components of $G\setminus (X\cup \{v\})$,
where every two vertices in $X$ are joined by an induced path of even length with 
interior in $V(A_1)$. 
\\
\\
(1) {\em For $1\le i\le k$ and for all distinct $x,x'\in X$, every induced path between $x,x'$ with interior in $A_i$ 
has a length that is even and at least four.}
\\
\\
For all distinct $x,x'\in X$, let $P_1(x,x')$ be an induced path of even length with
interior in $V(A_1)$. It follows that $P_1(x,x')$ has length at least four, since $v$ is adjacent to $x,x'$ and $G$
has girth at least five. For $2\le i\le k$, if $x,x'$ both have a neighbour in $A_i$, let $P_i(x,x')$
be an induced path between $x,x'$ with interior in $A_i$. Thus $P_i(x,x')$ has length at least three, for the same reason.
Since $P_1(x,x')\cup P_i(x,x')$ is an induced cycle of length at least seven, it follows that $P_i(x,x')$ has even length, 
and length at least four, for all choices of $x,x'$ that both have a neighbour in $A_i$; and therefore every 
induced path between $x,x'$ with interior in $A_i$ has even length at least four. 
From the minimality of $X$,
it follows that every vertex in $X$ has a neighbour in $A_i$ for $2\le i\le k$; and so, by the same argument with $A_1,A_2$
exchanged, every induced path between $x,x'$ with interior in $A_1$ has even length at least four. This proves (1).

\bigskip

For $1\le i\le k$, let $G_i=G[V(A_i)\cup X\cup \{v\}]$. 
\\
\\
(2) {\em For $1\le i\le k$, there is a three-colouring  
$\phi_i:V(G_i)\rightarrow \{1,2,3\}$ with $\phi_i(v)=1$ and $\phi_i(x)=2$ for all $x\in X$.}
\\
\\
From the inductive hypothesis, $G_i$ admits a three-colouring 
$\phi_i:V(G_i)\rightarrow \{1,2,3\}$ with $\phi_i(v)=1$. Choose $\phi_i$ such that $\phi_i(x)=3$ for as few vertices $x\in X$
as possible. We claim that $\phi_i(x)=2$ for all $x\in X$. To see this, let $X_2$ be the set of $x\in X$
with $\phi_i(x)=2$, and let $X_3$ be the set of $x\in X$
with $\phi_i(x)=3$. Thus $X_2\cup X_3=X$. Let $H$ be the subgraph of $G_i$ induced on the set of vertices $u\in V(G_i)$
with $\phi_i(u)\in \{2,3\}$.  Suppose that $X_3\ne \emptyset$, and let $C$ be a component of $H$ that contains a vertex of $X_3$.
By exchanging colours in $C$, the choice of $\phi_i$ implies that some vertex of $X_2$ belong to $C$, and so 
there is a minimal induced path $P$
of $H$ between $X_2,X_3$, which therefore has odd length; but from the minimality of $P$, all internal vertices of $P$ 
belong to $A_i$, contradicting (1). This proves (2).

\bigskip

From (2) it follows that $G$ is three-colourable. This proves \ref{decompthm}.~\bbox

\section{Pentagraphs that contain large parts of the Petersen graph}

In this section we prove part of \ref{decompthm}. A {\em clique cutset} of $G$ is a clique $X$ of $G$ such that $G\setminus X$ is 
disconnected. In a pentagraph $G$, every clique has cardinality at most two, and so if $G$ has an edge and $G$ admits a 
clique cutset, then $G$ admits a strong parity star-cutset. If $P$ is an induced path, we denote the set of internal vertices of 
$P$ by $P^*$. Two disjoint subsets $X,Y$ of $G$ are {\em anticomplete} if there are no edges between $X,Y$. 
Let us say two nonadjacent vertices $s,t$ of a graph $H$ are {\em linked} if there are induced paths $Q_1,Q_2$
of $H$ both with ends $s,t$ and both of length at least three, with lengths of different parity. We say that $s,t$ are 
{\em odd-linked} if there is an induced path 
of $H$ with ends $s,t$ and with odd length at least five.

We begin with:
\begin{thm}\label{growing}
Let $G$ be a pentagraph that does not admit a clique cutset, and let $H$ be an induced subgraph of $G$ with $H\ne G$
and $|V(H)|\ge 3$. 
Then either
\begin{itemize}
\item there is a vertex $v\in V(G)\setminus V(G)$ with at least three neighbours in $V(H)$ (and therefore every two 
neighbours of $v$ in $V(H)$
have distance at least three in $H$); or 
\item there exist  nonadjacent $s,t\in V(H)$, and a vertex $v\in V(G)\setminus V(G)$ adjacent to $s,t$
and with no other neighbours in $V(H)$ (and therefore $s,t$ have distance at least three in $H$ and are not odd-linked in $H$); or
\item no vertex in $V(G)\setminus V(H)$ has more than one neighbour in $V(H)$, and 
there exist nonadjacent $s,t\in V(H)$, not linked in $H$, and an induced path $P$ of $G$ 
with length at least three, with ends $s,t$ and with $P^*\subseteq V(G)\setminus V(H)$, such that 
every vertex of $H$ with a neighbour in $P^*$ is adjacent to both $s,t$.
\end{itemize}
\end{thm}
\Proof
If some $v\in V(G)\setminus V(H)$ has at least three neighbours in $V(H)$ then 
the first bullet is satisfied, since $G$ has girth at least five. If some $v\in V(G)\setminus V(H)$ has exactly two 
neighbours $s,t$ 
in $V(H)$, and $Q$ is an induced path
of $H$ with ends $s,t$ and with odd length at least five, then  adding $v$ to $Q$
makes a long odd hole $G$, a contradiction; so $s,t$ are not odd-linked and the 
second bullet is satisfied. Thus we may assume that each vertex in $V(G)\setminus V(H)$ has at most one neighbour in $V(H)$.

Let $C$ be a component of $G\setminus V(H)$. Since $G$ does not admit a clique cutset, and $|V(H)|\ge 3$, it follows that
there exist nonadjacent vertices in $V(H)$ both with neighbours in $V(C)$. Thus there is an induced path $P$ with 
$P^*\subseteq V(C)$, with ends nonadjacent vertices of $C$. Choose $P$ with $P^*$ minimal, and let its ends be $s,t$. 
Since no vertex in $V(G)\setminus V(H)$ has more than one neighbour in $V(H)$, it follows that 
$P$ has length at least three.
If $v\in V(H)$ has a neighbour in $P^*$, and $v$ is nonadjacent to $s$ say, then from the minimality of $P^*$, it follows
that $v$ has only one neighbour in $P^*$ and that neighbour is adjacent to $t$. Hence $v$ is nonadjacent to $t$, since
$G$ has girth at least five; and this contradicts the minimality of $P^*$. So every vertex of $H$ with a neighbour in $P^*$ 
is adjacent to both $s,t$. Suppose that $s,t$ are linked in $H$, and so 
there are induced paths $Q_1, Q_2$ of $H$ between $s,t$, both of length at least three, and with lengths of 
different parity. Thus neither contains a vertex
adjacent to both $s,t$, and so $Q_1^*, Q_2^*$ are both anticomplete to $P^*$. Consequently both $P\cup Q_1, P\cup Q_2$
are long holes, and one has odd length, a contradiction; and so $s,t$ are not linked, and the 
third bullet is satisfied. This proves
\ref{growing}.~\bbox

We deduce:
\begin{thm}\label{Petersen}
Let $G$ be a pentagraph that has an induced subgraph isomorphic to the Petersen graph. Then either $G$ is isomorphic to 
the Petersen graph, or $G$ admits a clique cutset.
\end{thm}
\Proof
Let $H$ be an induced subgraph of $G$ isomorphic to the Petersen graph. No two vertices of $H$ have distance at least three 
in $H$. Moreover, every two nonadjacent vertices of $H$ are linked in $H$.
The result follows from \ref{growing}. This proves \ref{Petersen}.~\bbox

\begin{figure}[H]
\centering
\begin{tikzpicture}[scale=1,auto=left]
\tikzstyle{every node}=[inner sep=1.5pt, fill=black,circle,draw]
\def\r{2}
\def\s{1}
\node (v1) at ({\r*cos(90)}, {\r*sin(90)}) {};
\node (v2) at ({\r*cos(162)}, {\r*sin(162)}) {};
\node (v3) at ({\r*cos(234)}, {\r*sin(234)}) {};
\node (v4) at ({\r*cos(306)}, {\r*sin(306)}) {};
\node (v5) at ({\r*cos(18)}, {\r*sin(18)}) {};
\node (v6) at ({\s*cos(90)}, {\s*sin(90)}) {};
\node (v7) at ({\s*cos(162)}, {\s*sin(162)}) {};
\node (v8) at ({\s*cos(234)}, {\s*sin(234)}) {};
\node (v9) at ({\s*cos(306)}, {\s*sin(306)}) {};
\node (v10) at ({\s*cos(18)}, {\s*sin(18)}) {};

\foreach \from/\to in {v1/v2,v2/v3,v3/v4,v4/v5,v5/v1,v6/v8,v7/v9,v8/v10,v9/v6,v10/v7,v1/v6,v2/v7,v3/v8,v4/v9,v5/v10}
\draw [-] (\from) -- (\to);

\begin{scope}[shift ={(7,0)}]
\def\r{1.6}
\node (v1) at ({-\r}, {\r}) {};
\node (v2) at ({0}, {\r}) {};
\node (v3) at ({\r}, {\r}) {};
\node (v4) at ({\r}, {0}) {};
\node (v5) at ({\r}, {-\r}) {};
\node (v6) at ({0}, {-\r}) {};
\node (v7) at ({-\r}, {-\r}) {};
\node (v8) at ({-\r}, {0}) {};

\node (v9) at ({\s*cos(135)}, {\s*sin(135)}) {};
\node (v10) at ({\s*cos(45)}, {\s*sin(45)}) {};

\foreach \from/\to in {v1/v2,v2/v3,v3/v4,v4/v5,v5/v6,v6/v7,v7/v8,v8/v1,v2/v6,v4/v8,v1/v9,v9/v5,v3/v10,v10/v7}
\draw [-] (\from) -- (\to);

\tikzstyle{every node}=[]
\draw (v1) node [above]           {$1$};
\draw (v2) node [above]           {$2$};
\draw (v3) node [above]           {$3$};
\draw (v4) node [right]           {$4$};
\draw (v5) node [below]           {$5$};
\draw (v6) node [below]           {$6$};
\draw (v7) node [below]           {$7$};
\draw (v8) node [left]           {$8$};
\draw (v9) node [left]           {$9$};
\draw (v10) node [right]           {$10$};
\end{scope}

\end{tikzpicture}

\caption{$\mathcal{P}$ and $\mathcal{P}^0$.} \label{fig:Petersen}
\end{figure}
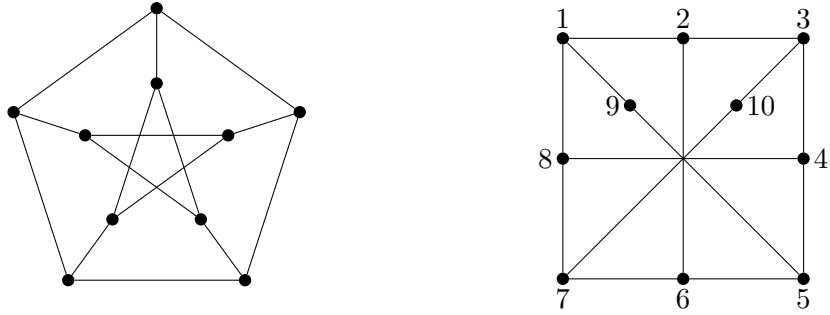

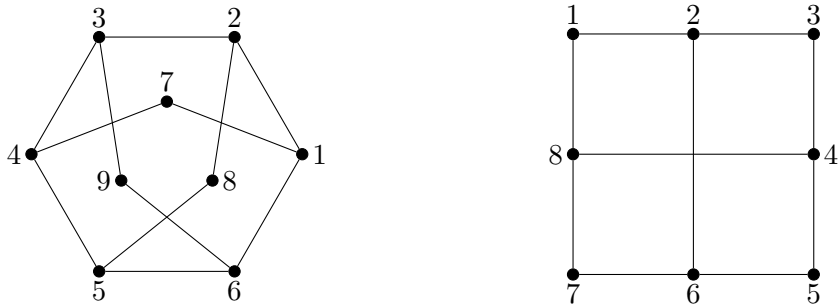
\begin{figure}[H]
\centering
\begin{tikzpicture}[scale=1,auto=left]
\tikzstyle{every node}=[inner sep=1.5pt, fill=black,circle,draw]

\def\r{1.8}
\def\s{.7}
\node (v1) at ({\r*cos(0)}, {\r*sin(0)}) {};
\node (v2) at ({\r*cos(60)}, {\r*sin(60)}) {};
\node (v3) at ({\r*cos(120)}, {\r*sin(120)}) {};
\node (v4) at ({\r*cos(180)}, {\r*sin(180)}) {};
\node (v5) at ({\r*cos(240)}, {\r*sin(240)}) {};
\node (v6) at ({\r*cos(300)}, {\r*sin(300)}) {};
\node (v7) at ({\s*cos(90)}, {\s*sin(90)}) {};
\node (v8) at ({\s*cos(330)}, {\s*sin(330)}) {};
\node (v9) at ({\s*cos(210)}, {\s*sin(210)}) {};

\foreach \from/\to in {v1/v2,v2/v3,v3/v4,v4/v5,v5/v6,v6/v1,v7/v1,v7/v4,v8/v2, v8/v5,v9/v3, v9/v6}
\draw [-] (\from) -- (\to);

\tikzstyle{every node}=[]
\draw (v1) node [right]           {$1$};
\draw (v2) node [above]           {$2$};
\draw (v3) node [above]           {$3$};
\draw (v4) node [left]           {$4$};
\draw (v5) node [below]           {$5$};
\draw (v6) node [below]           {$6$};
\draw (v7) node [above]           {$7$};
\draw (v8) node [right]           {$8$};
\draw (v9) node [left]           {$9$};

\begin{scope}[shift ={(7,0)}]
\tikzstyle{every node}=[inner sep=1.5pt, fill=black,circle,draw]
\def\r{1.6}
\node (v1) at ({-\r}, {\r}) {};
\node (v2) at ({0}, {\r}) {};
\node (v3) at ({\r}, {\r}) {};
\node (v4) at ({\r}, {0}) {};
\node (v5) at ({\r}, {-\r}) {};
\node (v6) at ({0}, {-\r}) {};
\node (v7) at ({-\r}, {-\r}) {};
\node (v8) at ({-\r}, {0}) {};

\foreach \from/\to in {v1/v2,v2/v3,v3/v4,v4/v5,v5/v6,v6/v7,v7/v8,v8/v1,v2/v6,v4/v8}
\draw [-] (\from) -- (\to);

\tikzstyle{every node}=[]
\draw (v1) node [above]           {$1$};
\draw (v2) node [above]           {$2$};
\draw (v3) node [above]           {$3$};
\draw (v4) node [right]           {$4$};
\draw (v5) node [below]           {$5$};
\draw (v6) node [below]           {$6$};
\draw (v7) node [below]           {$7$};
\draw (v8) node [left]           {$8$};
\end{scope}

\end{tikzpicture}

\caption{$\mathcal{P}^1$ and $\mathcal{P}^2$.} \label{fig:Petersen-1}
\end{figure}

Let $\mathcal{P}$ denote the Petersen graph, and let $\mathcal{P}^0,\mathcal{P}^1,\mathcal{P}^2$  denote the graphs obtained 
from $\mathcal{P}$
by deleting one edge,
one vertex, and two adjacent vertices respectively. (See figures \ref{fig:Petersen} and \ref{fig:Petersen-1}.)

\begin{thm}\label{P-0}
Let $G$ be a pentagraph that has an induced subgraph isomorphic to $\mathcal{P}^0$. Then either $G$ is isomorphic to 
$\mathcal{P}^0$, or 
$G$ admits a clique cutset.
\end{thm}
\Proof Let $H$ be an induced subgraph of $G$ isomorphic to $\mathcal{P}^0$, numbered as in figure \ref{fig:Petersen}. 
$(9,10)$ is the only pair of vertices of $H$ that have distance more than two in $H$.
Consequently no three vertices of $H$ pairwise have distance at least three; and every two vertices of $H$ with
distance at least three in $H$ 
are odd-linked in $H$. Moreover, every two nonadjacent vertices are linked in $H$, and so 
the result follows from \ref{growing}. This proves \ref{P-0}.~\bbox

\begin{thm}\label{P-1}
Let $G$ be a pentagraph that has an induced subgraph isomorphic to $\mathcal{P}^1$. Then either $G$ is isomorphic to
one of $\mathcal{P}, \mathcal{P}^0, \mathcal{P}^1$, or
$G$ admits a clique cutset.
\end{thm}
\Proof
Let $H$ be an induced subgraph of $G$ isomorphic to $\mathcal{P}^1$, numbered as in figure \ref{fig:Petersen-1}.
The only pairs of vertices that have distance at least three in $H$ are the pairs of vertices in $\{7,8,9\}$. Thus if the
first bullet of \ref{growing} holds then $G$ contains $\mathcal{P}$ and the result follows from \ref{Petersen}.
If the second bullet of \ref{growing} holds, then $G$ contains $\mathcal{P}_0$ and the result follows from \ref{P-0}.
Every two nonadjacent vertices of $H$ are linked, so the third bullet of \ref{growing} does not hold.
This proves \ref{P-1}.~\bbox

\begin{thm}\label{P-2}
Let $G$ be a pentagraph that has an induced subgraph isomorphic to $\mathcal{P}^2$. Then either $G$ is isomorphic to
one of $\mathcal{P}, \mathcal{P}^0, \mathcal{P}^1$, $\mathcal{P}^2$, or
$G$ admits a $P_3$-cutset or a strong parity star-cutset.
\end{thm}
\Proof
Let $H$ be an induced subgraph of $G$ isomorphic to $\mathcal{P}^2$, numbered as in figure \ref{fig:Petersen-1}.
The only pairs of vertices that have distance at least three in $H$ are $(1,5)$ and $(3,7)$, so the first bullet of 
\ref{growing} does not hold, and if the second bullet of \ref{growing} holds then $G$ contains $\mathcal{P}^1$ and the result follows from 
\ref{P-1} (since if $G$ has a clique cutset then it has a strong parity star-cutset). Thus we may assume that 
the third bullet of \ref{growing} holds, and in particular, 
no vertex in $V(G)\setminus V(H)$ has more than one neighbour in $V(H)$.

Let us say the four sets 
$$\{1,2,3\}, \{3,4,5\}, \{5,6,7\}, \{7,8,1\}$$
are the {\em sides} of $H$.
We may assume that $H\ne G$, and $G$ does not admit a $P_3$-cutset, and so there is a connected subgraph $F$ of $G\setminus V(H)$
such that $N(F)$ is not a subset of any side of $H$, where $N(F)$ denotes the set of vertices in $H$ with a 
neighbour in $V(F)$. Choose $F$ with $|V(F)|$ minimal. Since $N(F)$ is not a clique, there is an induced path $P$ with ends $s,t\in N(F)$,
nonadjacent; and as in \ref{growing}, by choosing $P$ with $P^*$ minimal it follows that $s,t$ are not linked in $H$, and
no vertex of $H$ has a neighbour in $P^*$ except $s,t$ and possibly a common neighbour of $s$ and $t$. 
The only nonadjacent pairs of vertices of $H$
that are not linked are $(1,3), (3,5), (5,7), (7,1)$, so from the symmetry we may assume that
$s=1$ and $t=3$.
Thus no vertex of $V(H)$ has a neighbour in $P^*$ except $1,3$ and possibly $2$. Since $N(F)\not\subseteq \{1,2,3\}$,
there is an induced path $Q$ of $G$ with interior in $V(F)$
with one end in $P^*$ and the other in $\{4,5,6,7,8\}$. Thus $Q$ has length at least two, and $F=P^*\cup Q^*$, from the 
minimality of $|V(F)|$.
Let the vertices of $Q$ be $q_1\DD q_2\CC q_k$ in order, where $k\ge 3$, and $q_1\in P^*$,
and $q_k\in \{4,5,6,7,8\}$. From the symmetry we may assume that $q_k\in \{4,5,6\}$. The vertices 
$1,2,3$ may have neighbours in $Q^*$, but since no vertex
has more than one neighbour in $V(H)$, it follows from the minimality of $|V(F)|$ that $q_k$ is the only vertex in $\{4,5,6,7,8\}$
with a neighbour in $Q^*$.
\\
\\
(1) {\em $q_k=4$.}
\\
\\
There is an induced path $R$ between $1,q_k$ with interior in $P^*\cup Q^*$, and it has length at least three since
no vertex in $V(G)\setminus V(H)$ has more than one neighbour in $V(H)$.
If $q_k= 5$, then
one of
$$1\DD R\DD 5\DD  4\DD 8\DD 1,1\DD R\DD 5\DD 6\DD 7\DD 8\DD 1$$
is a long odd hole, a contradiction; and if
$q_k= 6$, then
one of
$$1\DD R\DD 6\DD 5\DD 4\DD 8\DD 1,1\DD R\DD 6\DD 7\DD 8\DD 1$$
is a long odd hole, a contradiction. Since $q_k\in \{4,5,6\}$, this proves (1).

\bigskip

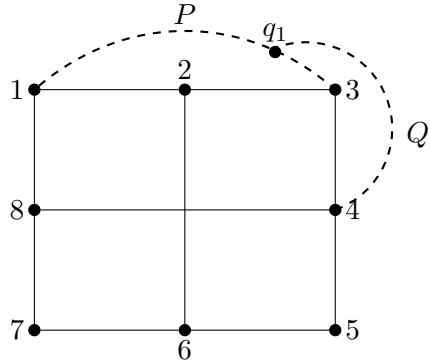
\begin{figure}[h!]
\centering
\begin{tikzpicture}[scale=1,auto=left]
\tikzstyle{every node}=[inner sep=1.5pt, fill=black,circle,draw]

\tikzstyle{every node}=[inner sep=1.5pt, fill=black,circle,draw]
\def\r{1.6}
\def\s{2}
\node (v1) at ({-\s}, {\r}) {};
\node (v2) at ({0}, {\r}) {};
\node (v3) at ({\s}, {\r}) {};
\node (v4) at ({\s}, {0}) {};
\node (v5) at ({\s}, {-\r}) {};
\node (v6) at ({0}, {-\r}) {};
\node (v7) at ({-\s}, {-\r}) {};
\node (v8) at ({-\s}, {0}) {};
\node (q1) at (1.2,2.1) {};

\foreach \from/\to in {v1/v2,v2/v3,v3/v4,v4/v5,v5/v6,v6/v7,v7/v8,v8/v1,v2/v6,v4/v8}
\draw [-] (\from) -- (\to);

\tikzstyle{every node}=[]
\draw (v1) node [left]           {$1$};
\draw (v2) node [above]           {$2$};
\draw (v3) node [right]           {$3$};
\draw (v4) node [right]           {$4$};
\draw (v5) node [right]           {$5$};
\draw (v6) node [below]           {$6$};
\draw (v7) node [left]           {$7$};
\draw (v8) node [left]           {$8$};
\draw (q1) node [above]           {$q_1$};

\draw[dashed, thick] (v1) to [bend left=40] (v3);
\draw[dashed,thick] (v4) arc (-70:110:1.15);
\tikzstyle{every node}=[]
\node at (0,2.6) {$P$};
\node at (3.1,1) {$Q$};
\end{tikzpicture}

\caption{Step (2) of the proof of \ref{P-2}.} \label{fig:P2jumps}
\end{figure}
\noindent
(2) {\em $2$ has no neighbour in $P^*\cup Q^*$, and $P$ has length three.}
\\
\\
(See figure \ref{fig:P2jumps}.) Suppose that $2$ has a neighbour in $P^*\cup Q^*$; then there is an induced path $R$ between $2, 4$ with interior in 
$P^*\cup Q^*$, and it has length at least three. 
But then one of 
$$2\DD R\DD 4\DD 5\DD 6\DD 2, 2\DD R\DD 4\DD 8\DD 7\DD 6\DD 2$$
is a long odd hole, a contradiction. Thus $2$ has no neighbour in $P^*\cup Q^*$. Consequently
$$1\DD P\DD 3\DD 2\DD 1, 1\DD P\DD 3\DD 4\DD 8\DD 1$$
are both holes of length at least five, and one has odd length; so $P$ has length three. 
This proves (2).

\begin{figure}[h!]
\centering
\begin{tikzpicture}[scale=1,auto=left]
\tikzstyle{every node}=[inner sep=1.5pt, fill=black,circle,draw]

\tikzstyle{every node}=[inner sep=1.5pt, fill=black,circle,draw]
\def\r{1.6}
\def\s{2}
\node (v1) at ({-\s}, {\r}) {};
\node (v2) at ({0}, {\r}) {};
\node (v3) at ({\s}, {\r}) {};
\node (v4) at ({\s}, {0}) {};
\node (v5) at ({\s}, {-\r}) {};
\node (v6) at ({0}, {-\r}) {};
\node (v7) at ({-\s}, {-\r}) {};
\node (v8) at ({-\s}, {0}) {};
\node (y) at (-1.5, 2.5) {};
\node (x) at (1.5, 2.5) {};

\foreach \from/\to in {v1/v2,v2/v3,v3/v4,v4/v5,v5/v6,v6/v7,v7/v8,v8/v1,v2/v6,v4/v8, v1/y,x/y,x/v3}
\draw [-] (\from) -- (\to);

\tikzstyle{every node}=[]
\draw (v1) node [left]           {$1$};
\draw (v2) node [above]           {$2$};
\draw (v3) node [right]           {$3$};
\draw (v4) node [right]           {$4$};
\draw (v5) node [right]           {$5$};
\draw (v6) node [below]           {$6$};
\draw (v7) node [left]           {$7$};
\draw (v8) node [left]           {$8$};
\draw (x) node [above]           {$x_1$};
\draw (y) node [above]           {$y_1$};

\draw[dashed,thick] (v4) arc (-70:95:1.3);
\tikzstyle{every node}=[]
\node at (3.1,1) {$Q$};
\end{tikzpicture}

\caption{Step (3) of the proof of \ref{P-2}.} \label{fig:P2jumps2}
\end{figure}

Let the vertices of $P$ be $1\DD y_1\DD x_1\DD 3$ in order.
\\
\\
(3) {\em $1$ has no neighbour in $V(Q)$, and so $q_1=x_1$.}
\\
\\
Suppose that $1$ has a neighbour in $V(Q)$; then there is an induced path $R$ between $1, 4$ with $R^*\subseteq V(Q)$.
Since $N(R^*)$ is not a subset of a side of $H$, the minimality of $|V(F)|$ implies that $R^*=V(F)$, and in particular
$x_1\in R^*\subseteq V(Q)$. Since $x_1\notin Q^*$, it follows that $q_1=x_1\in R^*$, and so $Q$ is a subpath of $R$,
contradicting that $1$ has a neighbour in $V(Q)$. This proves (3).

\bigskip

\begin{figure}[h!]
\centering
\begin{tikzpicture}[scale=1,auto=left]
\tikzstyle{every node}=[inner sep=1.5pt, fill=black,circle,draw]

\tikzstyle{every node}=[inner sep=1.5pt, fill=black,circle,draw]
\def\r{1.4}
\def\s{2}
\node (v1) at ({-\s}, {\r}) {};
\node (v2) at ({0}, {\r}) {};
\node (v3) at ({\s}, {\r}) {};
\node (v4) at ({\s}, {0}) {};
\node (v5) at ({\s}, {-\r}) {};
\node (v6) at ({0}, {-\r}) {};
\node (v7) at ({-\s}, {-\r}) {};
\node (v8) at ({-\s}, {0}) {};
\node (y3) at (-1.2, 2) {};
\node (x3) at (1.2, 2) {};
\node (y2) at (-1.2, 2.3) {};
\node (x2) at (1.2, 2.3) {};
\node (y1) at (-1.2, 2.6) {};
\node (x1) at (1.2, 2.6) {};

\foreach \from/\to in {v1/v2,v2/v3,v3/v4,v4/v5,v5/v6,v6/v7,v7/v8,v8/v1,v2/v6,v4/v8, v1/y1,x1/y1,x1/v3, v1/y2,y2/x2,x2/v3,
v1/y3,y3/x3,x3/v3}
\draw [-] (\from) -- (\to);

\tikzstyle{every node}=[]
\draw (v1) node [above]           {$1$};
\draw (v2) node [above]           {$2$};
\draw (v3) node [above]           {$3$};
\draw (v4) node [below right]           {$4$};
\draw (v5) node [below]           {$5$};
\draw (v6) node [below]           {$6$};
\draw (v7) node [below]           {$7$};
\draw (v8) node [left]           {$8$};
\draw (x1) node [above]           {$x_1$};
\draw (y1) node [above]           {$y_1$};
\node at (3.5,1.3) {$D$};
\node at (-1.2,3.3) {$Y$};
\node at (1.2,3.3) {$X$};

\draw (3.5,1.3) ellipse (0.5 and 1.5);

\draw (-1.2,2.6) ellipse (0.4 and 1);
\draw (1.2,2.6) ellipse (0.4 and 1);

\draw (x1) to (3.3, 2.4);
\draw (x1) to (3.3, 2.2);
\draw (v4) to (3.3, .5);
\draw (v4) to (3.3, .3);
\draw (v3) to (2.5,1.5);
\draw (v3) to (2.5,1.3);
\draw (x2) to (1.8,2.4);
\draw (x2) to (1.8,2.2);
\draw (x3) to (1.8,2.1);
\draw (x3) to (1.8,1.9);

\end{tikzpicture}

\caption{Handles in the proof of \ref{P-2}.} \label{fig:P2jumps3}
\end{figure}
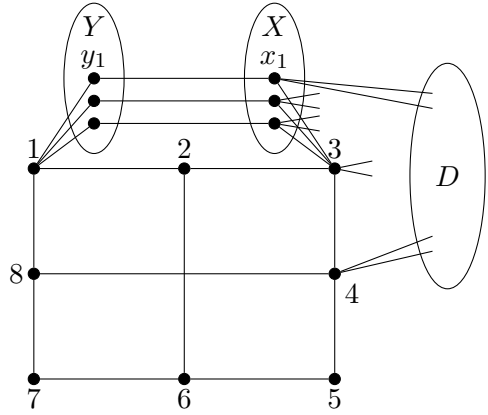

See figure \ref{fig:P2jumps2}. Let us say a {\em $1\DD 3$ handle} is an induced path $R$ of $G$ of length three between $1,3$ such that $R^*\cap V(H)=\emptyset$
and no vertex in $V(H)$ has a neighbour in $R^*$ except $1,3$. Thus $P$ is a $1\DD 3$ handle.
Let $X$ be the set of neighbours of $3$ that belong to $1\DD 3$ handles, and let $Y$ be the set of all neighbours of $1$
that belong to $1\DD 3$ handles.
Thus $x_1\in X$ and $y_1\in Y$, and $Y\cap V(Q)=\emptyset$ by (3). 
If some $y\in Y$ has a neighbour in $Q^*$, then $Q^*\cup \{y\}$ induces a connected subgraph of $G$ and $1,4$
both have neighbours in this subgraph, contrary to the minimality of $|V(F)|$. Thus $Y$ is anticomplete to $Q^*$, and therefore
$X\cap Q^*=\emptyset$.
Let $D$ be a connected induced subgraph of $G$, with $Q^*\subseteq V(D)$, maximal such that 
$V(D)\cap (X\cup \{3,4\})=\emptyset$ and no vertex in $Y\cup \{1,2,5,6,7,8\}$ has a neighbour in $V(D)$. It follows that $x_1,4$
both have a neighbour in $D$. (See figure \ref{fig:P2jumps3}.)

For every two vertices in $X\cup \{4\}$, there is a path between them of length four with middle vertex $1$; and this path is induced since
$G$ has girth at least five. We may assume that $X\cup \{3,4\}$ is not a strong parity star-cutset, and so $D$ is not a component
of $G\setminus (X\cup \{3,4\})$. Consequently there is a vertex $v\in V(G)\setminus V(D)$ with a neighbour in $V(D)$
and with $v\notin X\cup \{3,4\}$. From the maximality of $D$, it follows that $v$ has a neighbour in $Y\cup \{1,2,5,6,7,8\}$.
Since there is a path between this neighbour and $3$ of length at most three with vertex set in $V(H)\cup Y$, it follows that 
$3,v$ are nonadjacent. 
\\
\\
(4) {\em $v$ has a unique neighbour in $Y\cup \{1,2,5,6,7,8\}$.}
\\
\\
Suppose that $v$ has more than one neighbour in $Y\cup \{1,2,5,6,7,8\}$. Every vertex in $V(G)\setminus V(H)$ has at most
one neighbour in $V(H)$, as we saw earlier, so we may assume that $v$ is adjacent to some $y\in Y$. Choose $x\in X$ adjacent to $y$.
All neighbours of $v$ in $V(H)\cup X\cup Y$ pairwise have distance at least three in $G[V(H)\cup X\cup Y]$, and so
$v$ has no more neighbours in $Y$, and none in $\{x,1,2,3,8\}$. Since $v$ has two neighbours in 
$Y\cup \{1,2,5,6,7,8\}$, it has a unique neighbour in $\{5,6,7\}$ (say $u$), and therefore is nonadjacent to $4$.
The paths 
$$y\DD 1\DD 8\DD 7\DD 6\DD 5, y\DD 1\DD 8\DD 4\DD 5\DD 6, y\DD x\DD 3\DD 4\DD 8\DD 7$$
all have length five, and so $v$ cannot be adjacent to both ends of any of them, and so $u\ne 5,6,7$, a contradiction. 
This proves (4).

\bigskip

Let $u$ be the unique neighbour of $v$ in $Y\cup \{1,2,5,6,7,8\}$. 
Let $R$ be an induced path with interior in $V(D)\cup \{v\}$
between $x_1,u$, and
let $S$ be a minimal path with interior in $V(D)\cup \{v\}$ between $u$ and $\{3,4\}$. Thus one of $3,4$
is an end of $S$, and the other has no neighbour in $V(S)$. 
\begin{itemize}
\item 
If $u=5$, the union of $R$ with one of
$x_1\DD y_1\DD 1\DD 8\DD 7\DD 6\DD 5,x_1\DD y_1\DD 1\DD 2\DD 6\DD 5$
is a long odd hole. 
\item If $u=6$, the union of $R$ with one of
$x_1\DD y_1\DD 1\DD 8\DD 7\DD 6,x_1\DD y_1\DD 1\DD 2\DD 6$
is a long odd hole. 
\item If $u=7$, the union of $R$ with one of 
$x_1\DD y_1\DD 1\DD 8\DD 7,x_1\DD y_1\DD 1\DD 2\DD 6\DD 7$
is a long odd hole. 
\item If $u=8$, let $T$ be an induced path with interior in $V(D)\cup \{x_1, v\}$
between $3,u$; then $T$ has length at least three and the union of $T$ with one of 
$3\DD 2\DD 6\DD 7\DD 8, 3\DD 2\DD 1\DD 8$
is a long odd hole. 
\item If $u=2$, let $T$ be an induced path with interior in $V(D)\cup \{v\}$
between $4,u$; then $T$ has length at least three and the union of $T$ with one of
$4\DD 8\DD 7\DD 6\DD 2,4\DD 5\DD 6\DD 2$
is a long odd hole. 
\item If $u\in Y$ and $4$ is an end of $S$, then the union of $S$ with one of
$4\DD 5\DD 6\DD 2\DD 1\DD u, 4\DD 3\DD 2\DD 1\DD u$
is a long odd hole.
\item If $u\in Y$ 
and $3$ is an end of $S$, then $S$ has length at least three, and the union of $S$ with one of
$3\DD 2\DD 1\DD u, 3\DD 4\DD 8\DD 1\DD u$
is a long odd hole. 
\item If $u=1$ and $4$ is an end of $S$, then $S$ has length at least three, and the union of $S$ with one of
$4\DD 5\DD 6\DD 2\DD 1, 4\DD 3\DD 2\DD 1$
is a long odd hole.
\item If $u=1$ and $3$ is an end of $S$, then $S$ has length at least three, and the union of $S$ with one of
$3\DD 2\DD 1, 3\DD 4\DD 8\DD 1$
is an odd hole.
So $S$ is a path of length three between 1,3
with interior in $V(D)\cup \{v\}$, and no vertex of $H$ has a neighbour in $S^*$ except $1,3$. Consequently $S$ is a $1\DD 3$
handle, and so $V(D)\cup \{v\}$ contains a vertex in $X$, a contradiction.
\end{itemize}
Thus in all cases we obtain a contradiction. This proves \ref{P-2}.~\bbox

\section{Jumps across a pentagon}

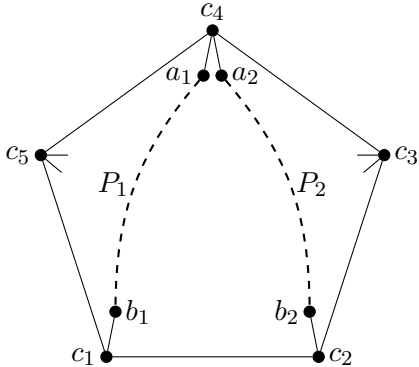
\begin{figure}[ht]
\centering
\begin{tikzpicture}[scale=1.2,auto=left]
\tikzstyle{every node}=[inner sep=1.5pt, fill=black,circle,draw]
\def\r{2}
\def\s{1}
\node (c4) at ({\r*cos(90)}, {\r*sin(90)}) {};
\node (c5) at ({\r*cos(162)}, {\r*sin(162)}) {};
\node (c1) at ({\r*cos(234)}, {\r*sin(234)}) {};
\node (c2) at ({\r*cos(306)}, {\r*sin(306)}) {};
\node (c3) at ({\r*cos(18)}, {\r*sin(18)}) {};
\node (a1) at ({\r*cos(90)-.1}, {\r*sin(90)-.5}) {};
\node (a2) at ({\r*cos(90)+.1}, {\r*sin(90)-.5}) {};
\node (b1) at ({\r*cos(234)+.1}, {\r*sin(234)+.5}) {};
\node (b2) at ({\r*cos(306)-.1}, {\r*sin(306)+.5}) {};

\foreach \from/\to in {c1/c2,c2/c3,c3/c4,c4/c5,c5/c1, a1/c4,a2/c4,b1/c1,b2/c2}
\draw [-] (\from) -- (\to);

\tikzstyle{every node}=[]
\draw (c2) node [right]           {$c_2$};
\draw (c3) node [right]           {$c_3$};
\draw (c4) node [above]           {$c_4$};
\draw (c5) node [left]           {$c_5$};
\draw (c1) node [left]           {$c_1$};
\draw (a1) node [left]           {$a_1$};
\draw (a2) node [right]           {$a_2$};
\draw (b1) node [right]           {$b_1$};
\draw (b2) node [left]           {$b_2$};

\draw[thick, dashed] (a1) to [bend right=20] (b1);
\draw[thick, dashed] (a2) to [bend left=20] (b2);

\node at (1.1,.3) {$P_2$};
\node at (-1.1,.3) {$P_1$};
\def\s{.3}
\draw (c3) to ({\r*cos(18)+\s*cos(180)}, {\r*sin(18)+\s*sin(180)});
\draw (c3) to ({\r*cos(18)+\s*cos(220)}, {\r*sin(18)+\s*sin(220)});
\draw (c5) to ({\r*cos(162)+\s*cos(0)}, {\r*sin(162)+\s*sin(0)});
\draw (c5) to ({\r*cos(162)+\s*cos(-40)}, {\r*sin(162)+\s*sin(-40)});

\end{tikzpicture}

\caption{Local jumps with a common end.} \label{fig:localjumps}
\end{figure}

In view of \ref{P-2}, we turn our attention to pentagraphs that do not contain $\mathcal{P}^2$ as an induced subgraph.
Let $G$ be a pentagraph, and let $C$ be a hole of length five in $G$. No vertex in $V(G)\setminus V(C)$ has more than one 
neighbour in $V(C)$, since $G$ has girth five.
Let $P$ be an induced path with both ends in $V(C)$, nonadjacent, and with no other vertices in $V(C)$. We call $P$ a 
{\em jump} over $C$. Let $P$ have ends $s,t$; then we call $P$ an {\em $s\DD t$ jump}, and if $c$ is the vertex of $C$ adjacent 
to both $s,t$, we say $P$ is a {\em jump across $c$}. 
If $P$ is an $s\DD t$ jump across $c$ and no vertex of $V(C)\setminus \{c,s,t\}$ has a neighbour in $P^*$, we say $P$ is a {\em local} jump.
If $P$ has length three we say that $P$ is a {\em short} jump. Then, clearly,
\begin{itemize}
\item all jumps have length at least three;
\item a jump $P$ is short if and only if no vertex of $V(C)\setminus V(P)$ has a neighbour in $P^*$;
\item short jumps are local;
\item local jumps have odd length.
\end{itemize}
We need to analyze which pairs of vertices of $C$ can be joined by short and local jumps.
We begin with:
\begin{thm}\label{localjumps}
Let $G$ be a pentagraph not containing $\mathcal{P}^2$, and let $C$ be a hole of length five in $G$.
If $P_1, P_2$ are local jumps over $C$ with exactly one common end $c$, 
then there is a short jump across $c$
with interior in $P_1^*\cup P_2^*$, and consequently neither of them is short.
\end{thm}
\Proof
(See figure \ref{fig:localjumps}.)
Suppose that there are local jumps $P_1,P_2$ over $C$ with exactly one common end $c$ say, and there is no short jump
across $c$ with interior in $P_1^*\cup P_2^*$. Choose such $P_1,P_2,c$ with $P_1^*\cup P_2^*$
minimal. (Note that $P_1^*\cap P_2^*$ may be nonempty.)
Let $C$ have vertices
$c_1\DD c_2\DD c_3\DD c_4\DD c_5\DD c_1$
in order, where $P_i$ is a $c_i\DD c_4$ jump for $i = 1,2$.
For $i = 1,2$, let $a_i, b_i$ be the vertices of $P_i$ adjacent to $c_i$ and $c_4$ respectively. For $i = 1,2$, let
$D_i=P_i^*\setminus \{a_i,b_i\}$. 
\\
\\
(1) {\em $D_1\cup \{b_1\}$ is disjoint from and anticomplete to $D_2\cup \{b_2\}$.}
\\
\\
Suppose not.
Since $b_2\notin V(P_1)$ (because $P_1$ is local) and vice versa, and $b_1,b_2$
are not adjacent, it follows that either there is a path of $G[D_1\cup D_2\cup \{b_1\}]$ 
from $b_1$ to $D_2$ or a path of $G[D_1\cup D_2\cup \{b_2\}]$ from $b_2$ to $D_1$,
and from the symmetry we may assume the first. Hence $D_2\ne \emptyset$, so $P_2$ is not short, and so $c_3$ has a 
neighbour in $D_2$. Consequently there is a path between $c_3$ and $\{c_1,c_5\}$ with interior in $D_1\cup D_2\cup \{b_1\}$. 
Let $Q$ be a minimal path from $c_3$ to one of $c_1,c_5$, with interior in $D_1\cup D_2\cup \{b_1\}$. It follows that
one of $c_1,c_5$ is an end of $Q$ and the other has no neighbour in $Q^*$. Moreover, neither of $c_2,c_4$ has a neighbour
in $Q^*$, and so $Q$ is a short jump. The choice of $P_1,P_2$ implies that $Q$ is not a short jump across $c_4$, and so $c_1$
is an end of $Q$ and $Q$ is a short jump across $c_2$.
There is no short jump across $c_1$ with interior in $P_1^*\cup Q^*$, since $c_2$ has no neighbour in $P_1^*\cup Q^*$;
and since $P_1^*\cup Q^*$ is a proper subset of $P_1^*\cup P_2^*$, this contradicts the minimality of $P_1^*\cup P_2^*$.
This proves (1).

\bigskip

If $a_1=a_2$, then since the paths $P_1,P_2$ have odd length (because they are local), (1) implies that $G[V(P_1\cup P_2)]$
is an odd hole, which therefore has length five; and so $P_1,P_2$ are both short. But then $G[V(C\cup P_1\cup P_2)]$
is isomorphic to $\mathcal{P}^2$, a contradiction. Thus $a_1\ne a_2$, and since $a_1\notin V(P_2)$ and vice versa,
(1) implies that $P_1^*, P_2^*$ are disjoint, and every edge between them has an end in $\{a_1,a_2\}$. 
Since $G[V(P_1\cup P_2)]$ is not a long odd hole, there is an edge between $P_1^*, P_2^*$; so from the symmetry we may assume that
$a_1$ has a neighbour in $D_2\cup \{b_2\}$, and therefore there is a $c_2\DD c_4$ jump  with interior in $\{a_1, b_2\}\cup D_2$.
Since $c_1,c_5$ have no neighbours in $\{a_1, b_2\}\cup D_2$, this $c_2\DD c_4$ jump is local, contrary to the minimality of 
$P_1^*\cup P_2^*$. This proves \ref{localjumps}.~\bbox

To complete the proof of \ref{decompthm} we need:
\begin{thm}\label{C5jumps}
Let $G$ be a pentagraph not containing $\mathcal{P}^2$, and let $C$ be a hole of length five in $G$. Then either:
\begin{itemize}
\item some vertex of $C$ has degree two; or
\item $G$ admits a $P_3$-cutset; or
\item $G$ admits a strong parity star-cutset.
\end{itemize}
\end{thm}
\Proof
We claim first that we may number the vertices of $C$ as $c_1\DD c_2\DD c_3\DD c_4\DD c_5\DD c_1$ in order, such that:
\\
\\
(1) {\em There are no short jumps across any of $c_3,c_4,c_5$, there are no local jumps across $c_4$, and
every local jump across $c_3$ or $c_5$ contains a vertex that is
in a short jump.}
\\
\\
Let $S$ be the set of $c\in V(C)$ such that there is a short jump across $c$, and let $L$ be the
set of $c\in V(C)$ such that there is a local jump across $c$. Thus $S\subseteq L$.
If $L$ is a clique, then (1) holds; so we may assume that $c_2,c_5\in L$, where
the vertices of $C$ are $c_1\DD c_2\DD c_3\DD c_4\DD c_5\DD c_1$
in order.
\ref{localjumps} implies that if $c,c'\in L$ are nonadjacent, then neither of them is in $S$; so $S$ is a clique, and every
vertex of $S$ is adjacent to every vertex of $L\setminus S$. In particular $S=\{c_1\}$, and $L= \{c_5,c_1,c_2\}$.
By \ref{localjumps} again, either every local jump across $c_2$
contains a vertex in a short jump across $c_1$, or every local jump across $c_5$ contains such a vertex; and from the symmetry
we may assume the second. This proves (1).

\begin{figure}[ht]
\centering
\begin{tikzpicture}[scale=1.2,auto=left]
\tikzstyle{every node}=[inner sep=1.5pt, fill=black,circle,draw]
\def\r{2}
\node (c4) at ({\r*cos(90)}, {\r*sin(90)}) {};
\node (c5) at ({\r*cos(162)}, {\r*sin(162)}) {};
\node (c1) at ({\r*cos(234)}, {\r*sin(234)}) {};
\node (c2) at ({\r*cos(306)}, {\r*sin(306)}) {};
\node (c3) at ({\r*cos(18)}, {\r*sin(18)}) {};
\def\s{.6}
\node (x3) at ({\r*cos(18)+\s*cos(171)}, {\r*sin(18)+\s*sin(171)}) {};
\node (y3) at ({\r*cos(18)+\s*cos(225)}, {\r*sin(18)+\s*sin(225)}) {};
\node (x5) at ({\r*cos(162)+\s*cos(9)}, {\r*sin(162)+\s*sin(9)}) {};
\node (y5) at ({\r*cos(162)+\s*cos(-45)}, {\r*sin(162)+\s*sin(-45)}) {};
\node (y1) at ({\r*cos(234)+\s*cos(27)}, {\r*sin(234)+\s*sin(27)}) {};
\node (x1) at ({\r*cos(234)+\s*cos(81)}, {\r*sin(234)+\s*sin(81)}) {};
\node (x2) at ({\r*cos(306)+\s*cos(99)}, {\r*sin(306)+\s*sin(99)}) {};
\node (y2) at ({\r*cos(306)+\s*cos(153)}, {\r*sin(306)+\s*sin(153)}) {};

\foreach \from/\to in {c1/c2,c2/c3,c3/c4,c4/c5,c5/c1, c1/x1,c1/y1,c2/x2,c2/y2,c3/x3,c3/y3,c5/x5,c5/y5,x1/x3,y1/y3,x2/x5,y2/y5}
\draw [-] (\from) -- (\to);
\draw[rotate=54] (-1.45,0) ellipse (8pt and 12pt);
\draw[rotate=126] (-1.45,0) ellipse (8pt and 12pt);
\draw[rotate=18] (1.45,0) ellipse (8pt and 12pt);
\draw[rotate=162] (1.45,0) ellipse (8pt and 12pt);

\tikzstyle{every node}=[]
\draw (c2) node [right]           {$c_2$};
\draw (c3) node [right]           {$c_3$};
\draw (c4) node [above]           {$c_4$};
\draw (c5) node [left]           {$c_5$};
\draw (c1) node [left]           {$c_1$};
\def\m{-.5}
\def\n{0}
\def\q{.4}
\draw[thick, dotted] (\m, \n) to ({\m + \q*cos(234)}, {\n+ \q*sin(234)});
\draw[thick, dotted] (-\m, \n) to ({-\m - \q*cos(234)}, {\n+ \q*sin(234)});

\def\t{1.45}
\node at ({\t*cos(234)}, {\t*sin(234)}){$X_1$};
\node at ({\t*cos(306)}, {\t*sin(306)}){$X_2$};
\node at ({\t*cos(18)}, {\t*sin(18)}){$X_3$};
\node at ({\t*cos(162)}, {\t*sin(162)}){$X_5$};
\draw (0,1) ellipse (.9 and 0.4);
\node at (0,1) {$D$};
\draw (c4)-- (-.2,1.2);
\draw (c4)-- (.2,1.2);

\end{tikzpicture}

\caption{The numbering of $C$.} \label{fig:numbering}
\end{figure}
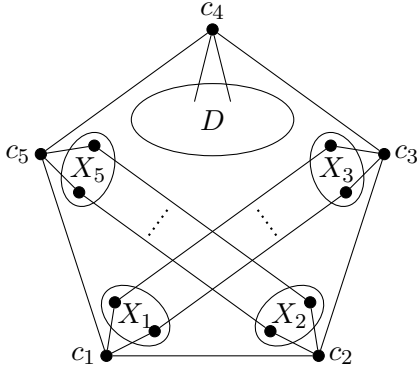

Let $X_1, X_3$ be the sets of vertices adjacent to $c_1,c_3$ respectively that are in short jumps across $c_2$ (thus, $X_1,X_3$
might be empty); and similarly let $X_2, X_5$ be the sets of vertices adjacent to $c_2,c_5$ respectively that are in short
jumps across $c_1$. (See figure \ref{fig:numbering}.) Let $X=X_1\cup X_2\cup X_3\cup X_5$. Thus $X$ is the set of all vertices that belong to the interior of short jumps.

Now, $c_4$ has no neighbour in $X\cup \{c_1,c_2\}$, but we may assume that
$c_4$ has degree at least three, and so there is a connected induced subgraph $D$ such that $c_4$ has a neighbour in $V(D)$ 
and $V(D)\cap (V(C)\cup X)=\emptyset$,
and $D$ is maximal with these properties. Let $N$ be the set of vertices in 
$V(C)\cup X$ that have a neighbour in $V(D)$; so $c_4\in N$.
\\
\\
(2) {\em $c_1,c_2\notin N$, and $N\cap (X_1\cup X_2)=\emptyset$.}
\\
\\
Suppose not; then from the symmetry we may assume that either $c_1$ or some member of $X_1$ belongs to $N$. Choose an induced
path $P$ between $c_4, c_1$ with interior in $D\cup X_1$.  
Since $P^* \setminus X_1 \subseteq D$, it follows that $P^* \cap X\subseteq X_1$; and so $|P^*\cap X|\le 1$, since $P$ is 
induced.
Suppose that $P$ is a local jump. From \ref{localjumps}, there is no short jump across $c_2$, and in particular $X_1=\emptyset$;
and so $V(P)\cap X=\emptyset$, contrary to (1).
Thus $P$ is not local. Let $Z$ be the set of vertices of $P$ that are not equal or adjacent to $c_1$ or to $c_4$. Thus 
one of $c_2,c_3$ has a neighbour in $Z$, since $P$ is not local and $G$ has girth five. If also $c_5$ has a neighbour in $P^*$, 
this neighbour also belongs to $Z$, and so there is a minimal path $Q$ from $c_5$ to one of $c_2,c_3$, with interior in $Z$.
Then no vertex of $C$ has a neighbour in $Q^*$ except the ends of $Q$, and so $Q$ is short, and therefore
two vertices of $Q^*$ belong to $X$; and since $Q^*\subseteq P^*$. But this contradicts that $Q^*\subseteq P^*$ and 
$|P^*\cap X|\le 1$.
So $c_5$ has no neighbour in $P^*$. If $c_2$ has a neighbour in $P^*$, choose a minimal path 
between $c_2,c_4$ with interior in $P^*$; then this is a local jump across $c_3$, containing no vertices in $X$, 
contrary to (1).
Thus $c_3$ has a neighbour 
in $P^*$, and hence in $Z$, and none of $c_2,c_4,c_5$ have a neighbour in $Z$; and so there is a short jump across $c_2$
with interior in $P^*$, and therefore two vertices of $P^*$ belong to $X$, a contradiction. This proves (2).

\bigskip

From (2) it follows that $N\subseteq X_3\cup X_5\cup \{c_3,c_4,c_5\}$. If $N\subseteq \{c_3,c_4,c_5\}$ then $G$ admits a 
$P_3$-cut, so we may assume from the symmetry that some $x_3\in X_3$ belongs to $N$. If also $c_5$ or some $x_5\in X_5$ has a neighbour in $D$,
then there is an induced path $Q$ between $c_3,c_5$ with interior in $V(D)\cup X_3\cup X_5$, and so neither of $c_1,c_2$
have a neighbour in it, and it is therefore a local jump across $c_4$, contrary to (1). Thus $N\subseteq X_3\cup \{c_3,c_4\}$.
But every two vertices in $X_3\cup \{c_4\}$ are joined by an induced path of length four with interior in 
$X_1\cup \{c_1,c_5\}$, and so $X_3\cup \{c_3,c_4\}$ is a strong parity star-cutset. This proves \ref{C5jumps}.~\bbox

Finally we deduce \ref{decompthm}, which we restate:
\begin{thm}\label{decompthm2}
Let $G$ be a pentagraph. Then either
\begin{itemize}
\item $G$ is bipartite; or
\item $G$ is isomorphic to the Petersen graph; or
\item $G$ has a vertex of degree at most two; or
\item $G$ admits a $P_3$-cutset or a strong parity star-cutset.
\end{itemize}
\end{thm}
\Proof Since $\mathcal{P}^0, \mathcal{P}^1, \mathcal{P}^2$ all have vertices of degree two, the result is true by \ref{P-2}
if $G$ contains an induced subgraph isomorphic to $\mathcal{P}^2$; so we assume it does not. We may assume that $G$ is not 
bipartite, and so it has a hole of length five. But then the result follows from \ref{C5jumps}. This proves \ref{decompthm2}.~\bbox

\section{Construction?}
Robertson's conjecture, that the Petersen graph is the only non-bipartite three-connected internally four-connected pentagraph, 
is false, but 
perhaps something like it is true. For instance, in~\cite{nelson} the following is shown:
\begin{thm}\label{cubic}
The Petersen graph is the only cubic non-bipartite three-connected pentagraph.
\end{thm}
A cubic three-connected pentagraph cannot admit a 
$P_3$-cutset (because the middle vertex of the $P_3$ must have a neighbour on either side of the cutset, by three-connectivity),
and cannot admit a strong parity star-cutset (because the ``strong'' condition implies that such a cutset would be a
$P_3$-cutset); and so in fact \ref{cubic} follows easily from \ref{decompthm}. 

Is there some hope of extending \ref{cubic} to larger classes of pentagraphs, or indeed to a construction for all pentagraphs?
In~\cite{plummer}, Plummer and Zha conjecture that every counterexample to Robertson's conjecture is ``close to bipartite'', and 
one might hope that this too would follow from \ref{decompthm}, but we do not see how to show it. The problem is, 
let $G$ be a pentagraph, with a parity star-cutset $X\cup\{x\}$, where $x$ is adjacent to every vertex in $X$, and
$A_1\LL A_k$ are the components of $G\setminus (X\cup \{x\})$. As in the derivation of \ref{mainthm} from \ref{decompthm},
we may assume that each vertex in $X$ has a neighbour in each of $A_1\LL A_k$, and every induced path with ends in $X$
and no other vertex has even length. Let $G_i$ be the subgraph induced on $V(A_i)\cup X\cup \{x\}$. We would like to apply
an inductive hypothesis that says each $G_i$ is ``close to bipartite'', but even if $G$ is three-connected and
internally four-connected, we do not know that $G_1\LL G_k$ are three-connected; for instance, they might have vertices of degree two.

On the other hand, a parity star-cutset is a ``reversible'' decomposition, that can be turned into something like a construction:
in the notation above, if we do not know that $G$ is a pentagraph, but we know that 
each of $G_1\LL G_k$ is a pentagraph, it folows that $G$ is indeed a pentagraph. So there is some hope here for a construction.

\section*{Acknowledgements}
This research was conducted during an Oberwolfach workshop on graph theory in January 2022, and the authors
are grateful to MFO for providing accommodations and facilities.
We would like to thank Xingxing Yu for introducing us to this problem (at an open problem session during the workshop) and for 
supplying background information.

\end{document}